\documentstyle[11pt,leqno]{article}
\newtheorem{theorem}{{\sc Theorem}}
\newcommand{\bt}{\begin{theorem}}
\newcommand{\et}{\end{theorem}}
\newcommand{\newsection}[1]{\setcounter{equation}{0} \setcounter{theorem}{0}
\section{#1}}

\newcommand{\NI}{\noindent}
\newcommand{\bea}{\begin{eqnarray}}
\newcommand{\eea}{\end{eqnarray}}

\def \spec#1 {\mathop{#1}}

\def \b #1 {\bf #1}
\newcommand {\nnb }{\nonumber}

\newcommand {\CC}{\centerline}

\newcommand{\ity}{\infty}
\newcommand{\raro}{\rightarrow}

\newcommand{\vsp}{\vskip 1em}

\newcommand{\be}{\begin{equation}}
\newcommand{\ee}{\end{equation}}
\newcommand{\ben}{\begin{eqnarray*}}
\newcommand{\een}{\end{eqnarray*}}

\begin{document}
\CC {\bf Singularity  for Bifractional and Trifractional  Brownian  Motions}
\CC{\bf Based on Their Hurst Indices}
\vsp
\CC{\bf B. L. S. Prakasa  Rao }
\CC{CR Rao Advanced Institute of Mathematics, Statistics}
\CC{ and Computer Science, Hyderabad, India} 
\vsp 
\NI{\bf Abstract:} We study sufficient conditions which ensure that the probability measures generated by two bifractional Brownian motions on an interval $[0,1]$ are singular with respect to each other and sufficient conditions for the probability measures generated by two trifractional Brownian motions on an interval $[0,1]$ are singular with respect to each other.
\vsp
\NI{\bf Keywords:} Singularity for measures; Bifractional Brownian motion; Trifractional Brownian motion; Hurst indices.
\vsp
\NI{\bf AMS Subject Classification:} 60G22; 60G30.
\newsection{Introduction}

A Gaussian stochastic process $\{X(t), -\ity < t <\ity\}$ with zero mean and covariance function 
$$ r(s,t)=\frac{1}{2}(|s|^{2H}+|t|^{2H}-|s-t|^{2H})$$
with $0<H<1$ is called a {\it fractional Brownian motion} with Hurst index $H.$
This process is a standard Brownian motion if $H=\frac{1}{2}.$ Probabilistic properties of fractional Brownian motion (fBm) and self-similar processes and their application to modeling phenomena with long-range dependence is discussed in Doukhan et al. (2003). Some recent  results in statistical inference for processes driven by fBm were surveyed in Prakasa Rao (2010). 

\newsection{Baxter type theorems}

\NI{\bf Theorem 2.1:} Suppose the process $\{X(t), a \leq t \leq b\}$ is a Gaussian process with mean function $m(t)$ and covariance function $r(s,t).$ Let $T=[a,b].$ Suppose the function $m(.)$ is differentiable and the function $r(s,t)$ is continuous on the square $T \times T$ with uniformly bounded (mixed) second derivatives just off the diagonal of the square $T \times T.$ Let
$$f_r(t) = D^{-}(t)-D^{+}(t)$$
where
$$D^{+}(t)= \lim_{s \raro t^{+}}\frac{r(t,t)-r(s,t)}{t-s}$$
and
$$D^{-}(t)= \lim_{s \raro t^{-}}\frac{r(t,t)-r(s,t)}{t-s}.$$
Then   
$$\lim_{n \raro \ity}\sum_{k=1}^{2^n}[X(a+\frac{k(b-a)}{2^n})-X(a+\frac{(k-1)(b-a)}{2^n})]^2=\int_a^bf_r(s)ds \;\;\mbox{a.s.}\;\;$$
\vskip10pt
The result stated above is due to Baxter (1956)  (cf. Rao (2000), p.249). As a special case of Theorem 2.1, we have the following result due to Levy.

\vskip10pt
\NI{\bf Theorem 2.2:} Consider a standard Brownian motion $\{W(t), 0 \leq t \leq T\}.$ Then
$$\lim_{n \raro \ity}\sum_{k=1}^{2^n}[W(\frac{kT}{2^n})-W(\frac{(k-1)T}{2^n})]^2= T \;\;\mbox{a.s.}\;\;$$
\vskip10pt

The following result is due to Decreusefond and Ustunel (1999) for a fractional Brownian motion (fBm).
\vskip10pt
\NI{\bf Theorem 2.3:} Let $\{W_H(t), 0 \leq t \leq T\}$ be a fractional Brownian motion with Hurst index $H \in (0,1)$ as defined above. Then, with probability one,
\begin{eqnarray*} 
\lim_{n \raro \ity}\sum_{k=1}^{2^n}|W_H(\frac{kT}{2^n})-W_H(\frac{(k-1)T}{2^n})|^p & = & 0\;\;\mbox{if}\;\; pH >1\\
& = & \ity \;\;\mbox{if}\;\; pH <1\\
& = & T \;\;\mbox{if}\;\; pH =1.
\end{eqnarray*}
\vskip10pt
Observe that by choosing $H=\frac{1}{2}$ and $p=2$ in Theorem 2.3, we obtain $pH=1$ and we obtain Theorem 2.2. However, if the Hurst index $H > \frac{1}{2},$ then the the quadratic variation of the process $W_H(.)$  (the case when p=2) is zero since $pH >1.$ If $H < \frac{1}{2},$ then  the quadratic variation of the process $W_H(.)$ is infinite since $p=2, H< \frac{1}{2}$ and hence $pH <1.$
\vskip10pt
In a recent paper on estimation of the Hurst index for fBm, Kurchenko (2003) derived Baxter type theorems for a fBm. 

Let $f: (a,b) \raro R$ be a function and let $k$ be a positive  integer.
Let $\Delta_h^{(k)}f(t)$ denote the increment of $k$-th order of the function $f$ on an interval $[t,t+h] \subset (a,b)$ as given, namely,
$$ \Delta_h^{(k)}f(t)= \sum_{i=0}^k (-1)^i {k \choose i}f(t+\frac{i}{k}h).$$
For any $m \geq 0$, positive integer $k \geq 1$ and $0<H<1,$ define

$$V_k(m,H)= \frac{1}{2}\sum_{i,j=0}^k(-1)^{i+j+1}{k\choose i}{k\choose j}|m+\frac{i-j}{k}|^{2H}.$$

It can be checked that $V_2(0,H)= 2^{2-2H}-1$ and
$$\Delta_1^{(2)}f(t)= f(t)-2f(t+\frac{1}{2})+f(t+1).$$
 
Kurchenko (2003) proved the following Baxter type theorem for second order increments for a fBm on the interval $[0,1].$ 
\vskip10pt
\NI{\bf Theorem 2.4 :} Let $\{W_H(t), t \in [0,1]\}$ be a fractional Brownian motion with Hurst index $H \in (0,1)$ as defined above. Then, with probability one,
$$\lim_{n \raro \ity}\frac{1}{n}\sum_{m=0}^{n-1}(\Delta_1^{(2)}W_H(m))^2= V_2(0,H)\;\;\mbox{a.s.}\;\;$$
\vskip10pt
In other words
$$\lim_{n \raro \ity}\frac{1}{n}\sum_{m=0}^{n-1}
[W_H(m)-2W_H(m+\frac{1}{2})+W_H(m+1)]^2 = V_2(0,H)\;\;\mbox{a.s.}\;\;$$
for any fBm on the interval $[0,1]$ with Hurst index $H \in (0,1).$
\vskip10pt

Gladyshev (1961) proved results similar to those of Baxter (1956) for a class of Gaussian processes in particular for a fractional Brownian motion. This result can also be used for estimating the Hurst index of a fractional Brownian motion. Gladyshev's theorem has been extended
by Klein and Gine (1975), Marcus and Rosen (1992), Kono (1969), and Norvaisa (2011). Vitasaari (2019) gives necessary and sufficient conditions for limit theorems for quadratic variations of Gaussian sequences. Houdre and Villa (2003) studied a class of Gaussian processes processes termed as bifractional Brownian motion. Ma (2013) introduced another class of Gaussian processes which are called  trifractional Brownian motion. We will now study singularity properties of these processes.

\newsection{Bifractional Brownian motion}

\NI{\bf Definition 3.1:} A {\it bifractional Brownian motion} $\{B_t^{H,K}, t \in [0,T]\}$ is a centered Gaussian process, starting from zero, with covariance function
$$R^{H,K}(t,s)= \frac{1}{2^K}((t^{2H}+s^{2H})^K-|t-s|^{2HK}), 0\leq s,t \leq T$$
where $0<H<1$ and $0<K\leq 1.$ If $K=1,$ then the process $B^{H,1}$ reduces to a fBm on the interval $[0,T],$ with Hurst index $H$ in the interval (0,1). 
\vsp
The process $B^{H,K}= \{B_t^{H,K}, t \in [0,T]\}$ is $HK$-self-similar. It is Holder continuous of order $\delta>0$ for any $0<\delta<HK.$ If $HK= \frac{1}{2},$ then the process is a finite quadratic variation process on any interval $[0,T]$ with the quadratic variation equal to a constant times $T$. It is a short memory process and is neither a Markov process nor a semimartingale. The process does not have stationary increments except in case $K=1.$ Furthermore for any $T>0,$
$$2^{-K}|t-s|^{2HK}\leq E(B_t^{H,K}-B_s^{H,K})^2\leq 2^{1-K}|t-s|^{2HK}, 0\leq s,t \leq T.$$
\vsp
For more on properties of a bifractional Brownian motion, see Russo and Tudor (2006), Tudor and Xiao (2007) and Lie and Nualart (2009). The bifractional Brownian motion was first introduced by Houdre and Villa (2003) for $0<H<1$ and $0<K\leq 1$ and later extended for $H\in (0,1), K\in (1,2)$ with $HK \in (0,1)$ in Bardina and Es-Sebaiy (2011).  Viitasaari (2019) proved the following Baxter type theorem for a bifractional Brownian motion.
\vsp
\NI{\bf Theorem 3.1:} Let $B^{H,K}=\{B^{H,K}_t, 0\leq t \leq t\}$ be a bifractional Brownian motion on the interval $[0,T]$ with indices $H\in (0,1), K\in (0,2)$ and $HK\in (0,1).$ Let $\{\pi_n, n\geq 1\}$ be any sequence of partitions $\pi_n:0=t_0^n<t_1^n,\dots<t_n^n=T$ of the interval $[0,T]$ satisfying
$$|\pi_n|= o(\frac{1}{(\log n)^\gamma})$$
where 
$$ \gamma= \max(\frac{1}{2-2HK},1) \;\;\mbox{if}\;\;K\in (0,1)$$
and
$$\gamma= \frac{1}{\min(1,2H)+1-2HK} \;\; \mbox{if}\;\;K\in (1,2)$$
(Here $|\pi_n|$ denotes the norm of the partition $\pi_n$). Then
$$\sum_{t_k^n\in \pi_n}\frac{(B^{H,K}_{t^n_k}-B^{H,K}_{t^n_{k-1}})^2}{|t^n_k-t^n_{k-1}|^{2HK-1}}\raro 2^{1-K}T$$
almost surely as $n \raro \ity.$

\newsection{Trifractional Brownian motion}

\NI{\bf Definition 4.1:} A process $Z_{H,K}= \{Z_{H,K}(t), t\in [0,T]\}$ with indices $H,K \in (0,1)$ is called a {\it trifractional Brownian motion} if it is a centered Gaussian process with the covariance function
$$C_{H,K}(s,t)= t^{2HK}+s^{2HK}-(t^{2H}+s^{2H})^K. 0\leq t,s \leq T.$$
\vsp
A trifractional Brownian motion $Z_{H,K}= \{Z_{H,K}(t), t\in [0,T]\}$ is a $HK$-self-similar process and has non-stationary increments. Ma (2013) termed the process defined above as the trifractional Brownian motion.  Since the covariance structures of a bifractional Brownian motion and trifractional Brownian motion are similar, a better name for the Gaussian process defined by the Definition 4.1 would be, for instance,  bifractional Brownian motion of the second kind. However, we will stick to the terminology introduced by  Ma (2013). 
\vsp
Consider a trifractional Brownian motion $X$ defined on the interval $[0,T].$ Han (2021) studied limit results for sums of the type
$$2^{\alpha n}\sum_{k=1}^{2^n}(X(\frac{k}{2^n})-X(\frac{k-1}{2^n}))^2$$
as $n \raro \ity$ where $X$ is a trifractional Brownian motion. The following result is due to Han (2021). 
\vsp
\NI{\bf Theorem 4.1 :} Let $X$ be trifractional Brownian motion on the interval $[0,T]$ with $H,K \in (0,1)$. Then the following result hold:\\
(i) If $HK<\frac{1}{2},$ then
\bea
\lim_{n\raro \ity}2^{\alpha n}\sum_{k=1}^{2^n}(X(\frac{kT}{2^n})-X(\frac{(k-1)T}{2^n}))^2 &=& 0 \;\;\; \mbox{if}\;\;\;\alpha < 2HK \\\nnb
&=& \ity \;\;\; \mbox{if}\;\; \;\alpha >2HK\\\nnb
\eea
almost surely and\\
(ii)If $HK\geq \frac{1}{2},$ then
\bea
\lim_{n\raro \ity}2^{\alpha n}\sum_{k=1}^{2^n}(X(\frac{kT}{2^n})-X(\frac{(k-1)T}{2^n}))^2 &=& 0 \;\;\; \mbox{if}\;\; \;\alpha < 1 \\\nnb
&=& \ity \;\;\; \mbox{if}\;\;\; \alpha >1\\\nnb
\eea
almost surely. 
\vsp
As a consequence of this result, Han (2021) derived the following corollary.\\
\vsp
\NI{\bf Corollary 4.2:} Suppose $X$ is a trifractional Brownian motion on the interval $[0,1]$ with $H,K \in (0,1).$ Further  suppose that $HK \leq \frac{1}{2}.$ Then
\be
\lim_{n\raro \ity}2^{HK n}\sum_{k=1}^{2^n}(X(\frac{k}{2^n})-X(\frac{(k-1)}{2^n}))^2 =HK 
\ee
almost surely.
\vsp
\newsection{ Singularity }

It is well known that if $P$ and $Q$ are probability measures generated by two Gaussian processes, then these measures are either equivalent or singular with respect to each other (cf. Feldman (1958), Hajek (1958)). For a proof, see Rao (2000), p. 226.
\vsp
Let $\{X^i(t), t\in[0,T]\}, i=1,2$ be two bifractional fractional Brownian motions with
Hurst indices $(H_i,K_i), i=1,2$ with $H_i$ and $K_i$ satisfying the conditions stated in Theorem 3.1 of Viitsaari (2019). Suppose that $K_1\neq K_2.$ Since the processes $X^i, i=1,2$ defined on the interval $[0,T]$ are Gaussian, it follows that the probability measures generated by these Gaussian processes on the interval $[0,T]$ are either equivalent or singular with respect to each other. We will now prove that they are singular with respect to each other  from the results of Feldman and Hajek  stated above,.
\vskip10pt
\NI{\bf Theorem 5.1:} Let $\{X^i(t), 0\leq t \leq T\}$ be two bifractional Brownian motions with Hurst indices $(H_i,K_i), i=1,2$ with $H_i$ and $K_i$ satisfying the conditions stated in Theorem 3.1. Let $\gamma_1,\gamma_2$ as defined in Theorem 3.1. Suppose further that $\gamma_1=\gamma_2$ and that $K_1\neq K_2.$ Let $P_i^T$ be the probability measure generated by the process $\{X^i(t), t \geq 0\}$ for $i=1,2$ on the interval $[0,T].$ Then the probability measures $P_1^T$ and $P_2^T$ are singular with respect to each other.\\
\vsp
\NI{\bf Proof :} Let $\gamma= \gamma_1=\gamma_2.$ Applying Theorem 3.1 for the sequence of partitions $\{\pi_n, n\geq 1\}$ such that $$|\pi_n|= o(\frac{1}{(\log n)^{\gamma}}) $$ for the processes $X^i, i=1,2$ over the interval $[0,T],$ we obtain that
$$\sum_{t_k^n\in \pi_n}\frac{(X^i_{t^n_k}-X^i_{t^n_{k-1}})^2}{|t^n_k-t^n_{k-1}|^{2H_iK_i-1}}\raro 2^{1-K_i}T, i=1,2$$
almost surely as $n \raro \ity.$ Since $K_1\neq K_2 $ and since the convergence stated above is almost sure convergence under the corresponding probability measures, it follows that the measures $P_1^T$ and $P_2^t$ are singular with respect to each other.
\vsp
Let $\{X^i(t), t \in [0,1]\}, i=1,2$ be two trifractional fractional Brownian motions with
Hurst indices $(H_i,K_i), i=1,2$ such that $H_iK_i\leq \frac{1}{2}, i=1,2.$ Suppose that $H_1K_1 \neq H_2K_2.$ From the result of  feldman and Hajek stated above, it follows that the probability measures generated by these Gaussian processes are either equivalent or singular with respect to each other. We will now prove that they are singular with respect to each other.
\vsp
\NI{\bf Theorem 5.2:} Let $\{X^i(t), 0\leq t \leq 1\}, i=1,2$ be two trifractional Brownian motions with Hurst indices $(H_i,K_i), i=1,2$ such that $H_iK_i\leq \frac{1}{2}, i=1,2.$ Suppose that $H_1K_1 \neq H_2K_2.$ Let $P_i$ be the probability measure generated by the process $\{X^i(t), 0\leq t \leq 1\}$ for $i=1,2.$ Then the probability measures $P_1$ and $P_2$ are singular with respect to each other.
\vsp
\NI{\bf Proof :} Applying Corollary 4.2, we obtain that
$$\lim_{n\raro \ity}2^{H_iK_i n}\sum_{k=1}^{2^n}(X^i(\frac{k}{2^n})-X^i(\frac{(k-1)}{2^n}))^2 =H_iK_i \;\;\mbox{a.s.}$$
for $i=1,2.$ Since $H_1K_1\neq H_2 K_2 $ and since the convergence stated above is almost sure convergence under the corresponding probability measures, it follows that the measures $P_1$ and $P_2$ are singular with respect to each other.
\vsp
\NI{\bf Remarks:} Similar approaches have been adopted to prove that two fractional Brownian motions are singular if the corresponding Hurst indices are different in Prakasa Rao (2008), to check whether two sub-fractional Brownian motions are singular in Prakasa Rao (2012) and to check whether two fractional psuedo-diffusion processes are singular in Prakasa Rao (2016). Perrin et al. 2001) defined the $n$-th order fractional Brownian motion $B_{H,n}= \{B_{H.n}(t), t \geq 0\}$ as a centered Gaussian process with the covariance function
$$C_{H,n}(s,t)=(-1)^n\frac{C_H^n}{2}[|t-s|^{2H}-\sum_{j=0}^{n-1}(-1)^j {2H \choose j} ((\frac{t}{s})^js^{2H}+(\frac{s}{t})^jt^{2H})]$$
where $n \geq 1, H\in (n-1,n)$ and $C_H^n= (\Gamma(2H+1))|\sin(\pi H)|)^{-1}.$ For any $n > 1,$ the $n$-th order fractional Brownian motions has the Hurst index  $H$ in the interval $(n-1,n)$ which is outside the interval (0,1) which is the range for the Hurst index for the case of a fBm. However the $n$-th order fractional Brownian motion is $H$-self-similar and it reduces to the standard fBm for $n=1.$ It would interesting to derive Baxter type theorems for the $n$-th order fractional Brownian motion.
\vsp
\NI{\bf Acknowledgment:} This work was supported under the ``INSA Senior Scientist" scheme of the Indian National Science Academy at the CR Rao Advanced Institute of Mathematics, statistics and Computer science, Hyderabad, India.
\vskip10pt
\NI{\bf References}
\begin{description}
\item Baxter, G. (1956) A strong limit theorem for Gaussian processes, {\it Proc. Amer. Math. Soc.}, {\bf 7}, 522-527.

\item Bardina, X. and Es-Sebaiy, K. (2011) An extension of bifractional Brownian motion, {\it Comm. on Stoch. Anal.}, {\bf 5}, 333-340.

\item Decreusefond, L.  and Ustunel, A.S. (1999) Stochastic analysis of the fractional Brownian motion, {\it Potential Anal.}, {\bf 10}, 177-214.

\item Doukhan, P., Oppenheim, G. and Taqqu, M.S. (2003) {\it Theory of Long-Range dependence}, Birkhauser, Boston.

\item Feldman, J. (1958) Equivalence and perpendicularity of Gaussian processes, {\it Pacific J. Math.}, {\bf 8}, 699-708, correction {\it ibid.}, {\bf 9}, 1295-1296.

\item Gladyshev, E.G. (1961) A new limit theorem for stochastic processes with Gaussian increments, {\it Theory Probab. and its Appl.}, {\bf 6}, 52-61.

\item Hajek, J. (1958) On a property of normal distribution of any stochastic processes, {\it Czech Math. J.}, {\bf 8}, 610-618.

\item Han, Xiyue. (2021) A Gladyshev theorem for trifactional Brownian motion and $n$-th order fractional Brownian motion, arXiv:2105.02385v! [math.PR] 6 May 2021. 

\item Houdre, C. and Villa, J. (2003) {\it Contemporary Mathematics}, Amer. Math. Soc., {\bf 336}, 195-201.

\item Kono, N. (1969) Oscillation of sample functions in stationary gaussian processes. {\it  Osaka J. Math.}, {\bf 6}, 1-12.

\item Kurchenko, O.O. (2003) A consistent estimate of the Hurst parameter for a fractional Brownian motion, {\it Theor. Probab. and Math. Statist.}, {\bf 67}, 97-106.

\item Lei,P. and Nualart, D. (2009) A decomposition of the bifractional Brownian motion and applications, {\it Statist. probab. Lett.}, {\bf 79}, 619-624. 

\item  Ma, Chungsheng. (2013) The Schoenberg-Levy kernel and relationships among fractional Brownian motion, trifactional Brownian motion, and others, {\it Theory Probab. and its Appl.}, {\bf 57}, 619-632.

\item Marcus, M.B. and Rosen, J. (1992) $p$-variation of the local times of symmetric stable processes and of Gaussian processes with stationary increments, {\it Ann. Probab.}, ?, 1685-1713.

\item Norvaisa, R. (2011) A complement to Gladyshev's theorem, {\it Lithuanian Math. J.}, {\bf 51}, 26-35.

\item Perrin, E., Harba, R., Berzin-Joseph, C., Iribarren, I. and Bonami, A. (2001) $n$-th order fractional Brownian motion and fractional Gaussian noises, {\it IEEE Trans. Signal Processing.}, {\bf 49}, 1049-1059.

\item Prakasa Rao, B.L.S. (2008) Singularity of fractional Brownian motions with different Hurst indices, {\it Stoch. Anal. Appl.}, {\bf 26}, 334-337.

\item Prakasa Rao, B.L.S. (2010) {\it Statistical Inference for Fractional Diffusion Processes}, Wiley, London. 

\item Prakasa Rao, B.L.S. (2012) Singularity of subfractional Brownian motions with different Hurst indices, {\it Stoch. Anal.  Appl.} , {\bf 30} (2012), 538-542.

\item Prakasa Rao, B.L.S. (2016) Conditions for singularity for measures generated by two fractional pseudo-diffusion processes, {\it  Stoch. Anal. and Appl.}, {\bf 34} , 183-192.

\item Rao, M.M. (2000) {\it Stochastic Processes: Inference Theory}, Kluwer, Dordrecht.

\item Russo, F. and Tudor, C. (2006) On bifractional Brownian motion, {\it Stoch. Proc. Appl.}, {\bf 116}, 830-856.

\item Tudor, C. and Xiao, Y. (2007) Sample path properties of bifractional Brownian motion, {\it Bernoulli}, {\bf 13}, 1023-1052. 

\item Viitasaari, L. (2019) Necessary and sufficient conditions  for limit theorems for quadratic variations of Gaussian sequences, {\it Probability Surveys}, {\bf 16}, 62-98.
\end{description}

\end{document}